\documentclass[a4paper,oneside,10pt]{article}

\usepackage[english]{babel}    
\usepackage[utf8]{inputenc}    
\usepackage[T1]{fontenc}       
\usepackage[reqno]{amsmath}    
\usepackage{amssymb,amsthm}    
\usepackage{mathtools}         
\usepackage{url}               
\usepackage{icomma}            
\usepackage{units}             
\usepackage{enumerate}         
\usepackage{array}             
\usepackage[usenames]{color}   
\usepackage{graphicx}          
\usepackage{caption}           
\usepackage{times}             
\usepackage{titlesec}          
\usepackage{soul}
\usepackage[ruled]{algorithm2e}
\usepackage{todonotes}         
\usepackage[bookmarks, colorlinks=true, linkcolor=black, anchorcolor=black,
citecolor=black, filecolor=black, menucolor=black, runcolor=black,
urlcolor=black, pdfencoding=unicode]{hyperref}  
\usepackage[
paper=a4paper,
top=2.5cm,
bottom=2.5cm,
left=2cm,
right=2cm
]{geometry}  

\newcommand{\R}{\ensuremath{\mathbb{R}}}
\renewcommand{\b}{\boldsymbol}

\newcommand{\lap}{\nabla^2}

\DeclareMathOperator*{\argmin}{arg\,min}
\DeclareMathOperator{\metric}{d}
\DeclareMathOperator{\sgn}{sgn}

\newenvironment{itemize*}{\vspace{-6pt}\begin{itemize}\setlength{\itemsep}{0pt}\setlength{\parskip}{2pt}}{\end{itemize}}
\newenvironment{enumerate*}{\vspace{-6pt}\begin{enumerate}\setlength{\itemsep}{0pt}\setlength{\parskip}{2pt}}{\end{enumerate}}
\newenvironment{description*}{\vspace{-6pt}\begin{description}\setlength{\itemsep}{0pt}\setlength{\parskip}{2pt}}{\end{description}}

\newcommand{\Title}{Discretized boundary surface reconstruction}
\newcommand{\Author}{Mitja Jančič}
\newcommand{\AuthorTwo}{Viktor Cvrtila}
\newcommand{\AuthorThree}{Gregor Kosec}

\title{\Title}
\author{\Author}
\date{\today}
\hypersetup{pdftitle={\Title}, pdfauthor={\Author}, pdfcreator={\Author},
pdfproducer={\Author}, pdfsubject={}, pdfkeywords={}}  

\usepackage{multicol}
\usepackage{nth}

\newenvironment{Figure}
{\par\medskip\noindent\minipage{\linewidth}}
{\endminipage\par\medskip}

\titleformat{\section}{\scshape\normalsize\centering}{\Roman{section}.}{1em}{}
\titleformat{\subsection}{\itshape\normalsize\centering}{\Alph{subsection}.}{1em}{}
\titleformat{\subsubsection}{\itshape\normalsize}{\arabic{subsubsection})}{1em}{}

\titlespacing\section{0pt}{2.82mm plus 2mm minus 0mm}{-0.71mm}
\titlespacing\subsection{5.5mm}{2.12mm}{-1.06mm}
\titlespacing\subsubsection{3.18mm}{0.72mm}{-2.12mm}
\titlespacing\paragraph{5.08mm}{0.0mm}{1ex}
\begin{document}


\begin{center}
  \fontsize{24}{28}\selectfont
  \Title \\[3ex]
  \fontsize{11}{11}\selectfont
  \Author$^{1, 2}$,  \AuthorTwo$^{3}$,  \AuthorThree$^{2}$ \\[0.71mm]
  \fontsize{10}{10}\selectfont
  $^1$ ``Jožef Stefan'' International Postgraduate School, Ljubljana,
  Slovenia \\[0.5mm]
  $^2$ ``Jožef Stefan'' Institute, Parallel and Distributed Systems
  Laboratory, Ljubljana, Slovenia \\[0.5mm]
  $^3$ Faculty of Mathematics and Physics, University of Ljubljana,
  Ljubljana, Slovenia \\[1mm]
  
  \href{mailto:Mitja.Jancic@ijs.si}{mitja.jancic@ijs.si},
  \href{mailto:cvrtilaviktor@gmail.com}{cvrtilaviktor@gmail.com},
  \href{mailto:Gregor.Kosec@ijs.si}{gregor.kosec@ijs.si}
\end{center}

\vspace{1ex}

\begin{multicols}{2}
  
  \fontsize{9}{10}\selectfont
  {\bfseries\noindent
    
    Abstract -- Domain discretization is an essential part of the solution procedure in numerical simulations. Meshless methods simplify the domain discretization to positioning of nodes in the interior and on the boundary of the domain. However, generally speaking, the shape of the boundary is often undefined and thus needs to be constructed before it can be discretized with a desired internodal spacing. Domain shape construction is far from trivial and is the main challenge of this paper. We tackle the simulation of moving boundary problems where the lack of domain shape information can introduce difficulties. We present a solution for 2D surface reconstruction from discretization points using cubic splines and thus providing a surface description anywhere in the domain. We also demonstrate the presented algorithm in a simulation of phase-change-like problem.
    }
  \fontsize{10}{11}\selectfont
  
  
  \section{Introduction}
  \label{sec:introduction}
  Tractable solutions to partial differential equations (PDEs) are not easily obtained. Often advanced mathematical procedures or a series of simplifications are needed to obtain a closed form solution to a real-life problem~\cite{plestenjak2010numericne}. Therefore, in practice, we often rely on numerical treatment that provides us with a numerical approximation. For that, different numerical methods for solving PDEs have been proposed. Most commonly used, e.g.\ Finite Difference Method~\cite{smith1985numerical}, Finite Element Method~\cite{zienkiewicz2000finite}, Finite Volume Method~\cite{eymard2000finite}, Boundary Element Method~\cite{aliabadi2002boundary}, require a mesh to operate, while meshless methods approximate the differential operators only using scattered nodes~\cite{belytschko1996meshless} as shown in Fig.~\ref{fig:example}. This is an important advantage as the node positioning is considered to be easier then mesh generation, however, far from trivial. For that reason, several dedicated node positioning algorithms emerged~\cite{li2000point, lohner2004general, slak2018generation, kosec2018local}. 
  
  Historically speaking, meshless methods were introduced in the 1990s. Since then, different numerical procedures have been proposed, e.g.\ meshless Element Free Galerkin~\cite{Belytschko1994}, the Local Petrov-Galerkin~\cite{atluri1998new}, h-p Cloud Method~\cite{duarte1996h} and others. In this paper, we will use the meshless generalization of the traditional finite difference method (FDM) -- the Radial Basis Function-generated Finite Differences (RBF-FD) originally proposed by Tolstykh~\cite{tolstykh2003using}. The RBF-FD has already been used in a vast variety of applications ranging from linear elasticity~\cite{slak2019refined}, 4-dimensional problems~\cite{janvcivc2021monomial}, geosciences~\cite{fornberg2015primer}, fluid mechanics~\cite{kosec2018local}, dynamic thermal rating of power lines~\cite{maksic2019cooling}, etc.
  
  The fact that the domain discretization in the context of meshless methods is heavily simplified makes the meshless methods very attractive in the context of moving boundary problems, e.g.\ phase change problems~\cite{ALVAREZHOSTOS2019104321}. Providing a good discretization of a moving phase front is no easy task, as satisfying the quasi-uniform internodal spacing $h$ is crucial to assure the stability of the numerical method~\cite{fornberg2015solving}. The only way to satisfy the quasi-uniform spacing $h$ condition on moving boundary problems is by repositioning the nodes from the domain. However, the node repositioning needs to bo performed with a minimum cost to domain shape distortion. To reduce the distortion, the domain boundary shape must be known even between the discretization points, meaning, a proper surface reconstruction algorithm from a set of boundary nodes is needed.
  
  Surface reconstruction has already been addressed in the context of numerical simulations -- Non-uniform rational basis splines (NURBS), often used in computer graphics for representing curves and surfaces~\cite{piegl1996nurbs}, are used in numerical simulations using the finite element analysis, e.g.\ Isogeometric analysis (IGA)~\cite{hughes2005isogeometric}. In this paper, we tackle surface reconstruction in two-dimensional domain space where the surface is represented as a two-dimensional curve. We provide an algorithm that uses cubic splines to reconstruct the domain shape from the discretization points and thus provides us with the domain shape information anywhere on the boundary. 

  \begin{Figure}
    \centering
    \includegraphics[width=\linewidth]{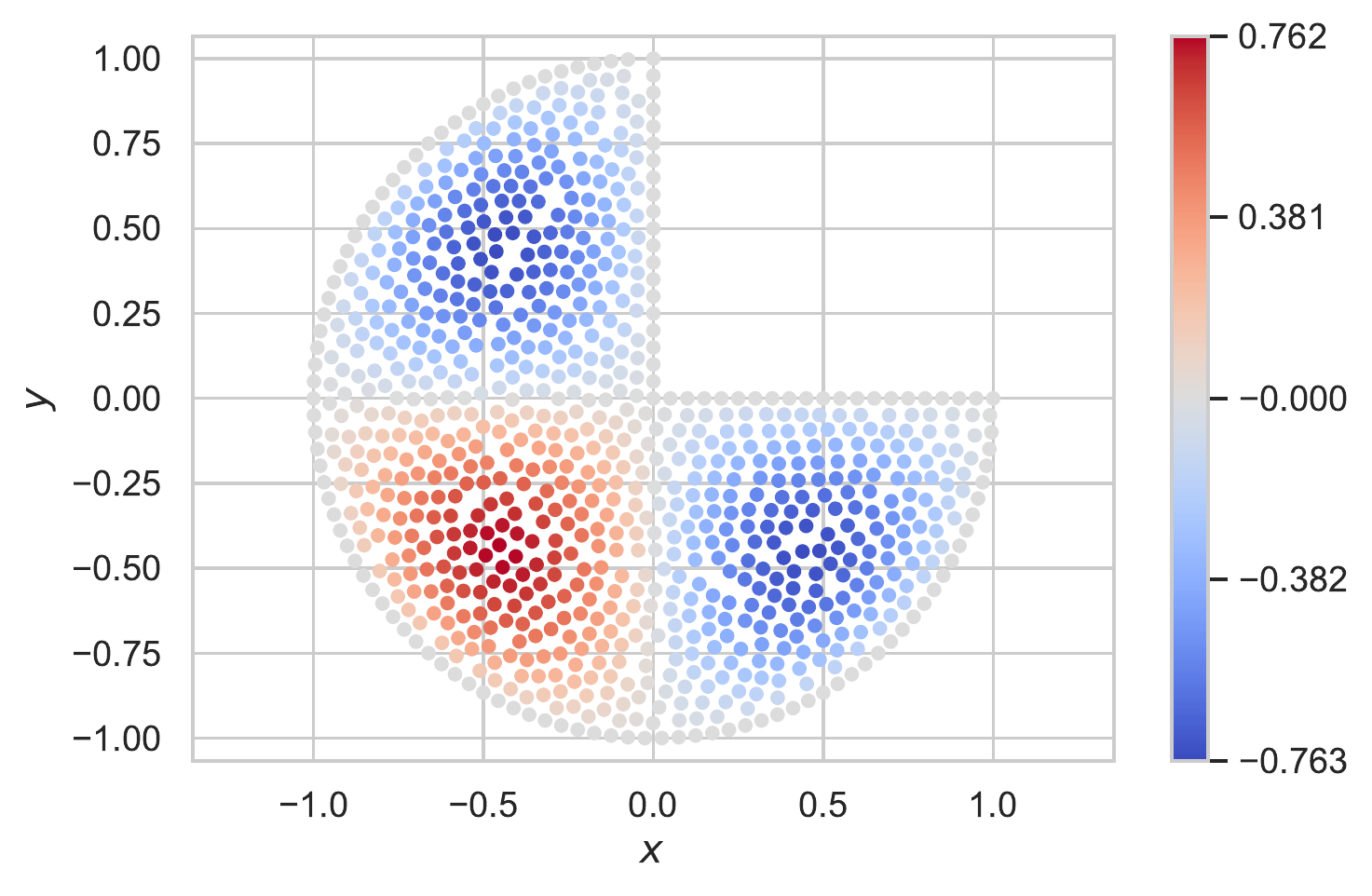}
    \captionof{figure}{Example solution of Poisson's problem with Dirichlet
      boundary conditions on $N=881$ scattered nodes in the domain.}
    \label{fig:example}
  \end{Figure}
  
  The paper is organized as follows: In section~\ref{sec:meshless} a short presentation to local strong form meshless methods, namely the RBF-FD is given. In section~\ref{sec:reconstruction} and section~\ref{sec:implementation} our proposed surface reconstruction algorithm is thoroughly explained. The algorithm is then used to solve a problem from section~\ref{sec:problem}, while results are presented in section~\ref{sec:results}. Final conclusions are gathered in the final section~\ref{sec:conclusions}.
  
  \section{Local strong form meshless methods}
  \label{sec:meshless}
  
  A general idea of the meshless methods is to use the local discretization points and construct an approximation of the considered field. This approximation is later used for manipulation with differential operators using the ansatz
  \begin{equation}
    \label{eq:ansatz}
    (\mathcal L u)(\b x) \approx \sum _{i}w_i^{\mathcal L} u(\b p_i).
  \end{equation}
  Here $\mathcal L$ is a differential operator and index $i$ runs over the set of closest neighbors $\b p _i$ of $\b x$. Equality of~\eqref{eq:ansatz} is enforced and the weights $\b w$ are computed. Different methods can be used to compute the weights, we will use the RBF-FD. Often used RBFs, e.g.\ Gaussians, include a shape parameter that can play a crucial role in the overall method stability~\cite{flyer2016role}. However, using Polyharmonic splines (PHS) and additionally augmenting them with polynomials helps overcome the stability issues~\cite{bayona2017role}.
  
  Some solution procedures, including the RBF-FD based, have been implemented using the object-oriented approach and C++’s strong template system. Node positioning, support selection, differential operator approximation, PDE discretization and other modules are all available as part of the Medusa library~\cite{medusa}, also used in this paper.
  
  \section{Surface reconstruction}
  \label{sec:reconstruction}
  Consider a moving boundary problem. The altering distance between the neighboring nodes can, firstly, increase to a point where the numerical methods become unstable due to the violation of quasi-uniform internodal spacing $h$ requirement, or secondly, become too large to achieve a desired accuracy of numerical solution in a specific domain area. To avoid such difficulties, repositioning the nodes in the domain is required, however, during the repositioning process the domain shape must be preserved as much as possible. Since the domain shape is, generally speaking, undefined, it first has to be constructed from a set of points before the boundary can be discretized. In this section we present a solution for 2D surface reconstruction from a set of points using cubic splines and provide complete information about the surface shape at hand.
  
  Suppose we are given a set of $k$ points
  $$ X = \{ \b x_i \in \R^2, \  i=0,\ldots, k-1 \}$$ representing the boundary $\partial \Omega$ of the domain $\Omega$, which is parametrized by a \emph{Jordan curve} $\b \gamma\colon [a, b] \to \R^2$. We do not have access to neither the domain nor the curve. Suppose there exist knots $\{t_j \in [a, b],\ j=0,\ldots, k-1\}$ such that $\b \gamma(t_j) = \b x_i$ for all $j$.
  The points are given in no particular order, i.e.~$0 \le i < j \le k-1$ does not necessarily imply that $t_i < t_j$.
  Our task is to find a curve $\tilde{\b \gamma}\colon [\tilde a, \tilde b] \to \R^2$ that approximates the original curve $\b \gamma$. 
  As there are many possible curves that interpolate $X$, obtaining the original curve is impossible without providing additional information or constraints.
  We, therefore, assume that the given points are \emph{dense enough} to adequately describe the curve in the following way.
  
  Let $Y = \b\gamma(\R) \subset \R^2$ be the image of the curve $\b \gamma$.
  Suppose $\b x_p$ and $\b x_q$ are neighboring points to $\b x_i$ in the sense that $t_i$ is the only knot between $t_p$ and $t_q$.
  The indices $p$ and $q$ are both dependant on the choice of $i$, i.e.\ $p=p(i)$ and $q=q(i)$, however, the explicit dependency is omitted in our writing.
  Let us then define an open neighborhood $U_i = \{\b \gamma(t),\ t_{p} < t < t_{q}\}$ of $\b x_i$ in $Y$. To help clarify the notation used, an illustration is provided in Fig.~\ref{fig:distance_illustration}.
  \begin{Figure}
    \centering
    \def\svgwidth{\columnwidth}
    \scalebox{0.8}{
\begingroup%
  \makeatletter%
  \providecommand\color[2][]{%
    \errmessage{(Inkscape) Color is used for the text in Inkscape, but the package 'color.sty' is not loaded}%
    \renewcommand\color[2][]{}%
  }%
  \providecommand\transparent[1]{%
    \errmessage{(Inkscape) Transparency is used (non-zero) for the text in Inkscape, but the package 'transparent.sty' is not loaded}%
    \renewcommand\transparent[1]{}%
  }%
  \providecommand\rotatebox[2]{#2}%
  \newcommand*\fsize{\dimexpr\f@size pt\relax}%
  \newcommand*\lineheight[1]{\fontsize{\fsize}{#1\fsize}\selectfont}%
  \ifx\svgwidth\undefined%
    \setlength{\unitlength}{226.77165354bp}%
    \ifx\svgscale\undefined%
      \relax%
    \else%
      \setlength{\unitlength}{\unitlength * \real{\svgscale}}%
    \fi%
  \else%
    \setlength{\unitlength}{\svgwidth}%
  \fi%
  \global\let\svgwidth\undefined%
  \global\let\svgscale\undefined%
  \makeatother%
  \begin{picture}(1,0.375)%
    \lineheight{1}%
    \setlength\tabcolsep{0pt}%
    \put(0,0){\includegraphics[width=\unitlength,page=1]{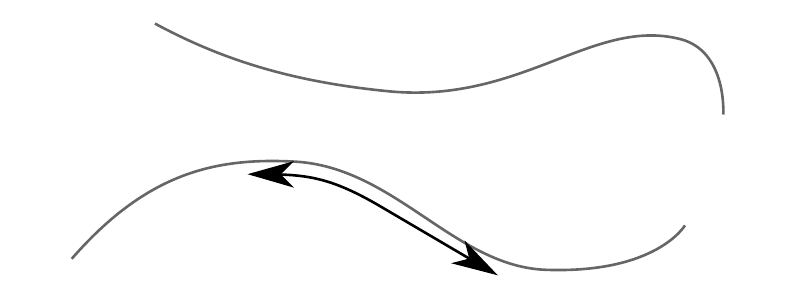}}%
    \put(0.44231067,0.0847913){\color[rgb]{0,0,0}\makebox(0,0)[lt]{\lineheight{1.25}\smash{\begin{tabular}[t]{l}$U_i$\end{tabular}}}}%
    \put(0.13376838,0.05185852){\makebox(0,0)[lt]{\lineheight{1.25}\smash{\begin{tabular}[t]{l}$\gamma$\end{tabular}}}}%
    \put(0,0){\includegraphics[width=\unitlength,page=2]{images/curve_distance.pdf}}%
    \put(0.51986251,0.1267524){\makebox(0,0)[lt]{\lineheight{1.25}\smash{\begin{tabular}[t]{l}$\b x_i$\end{tabular}}}}%
    \put(0.30790336,0.19325155){\makebox(0,0)[lt]{\lineheight{1.25}\smash{\begin{tabular}[t]{l}$\b x_{p(i)}$\end{tabular}}}}%
    \put(0.63276191,0.06593874){\makebox(0,0)[lt]{\lineheight{1.25}\smash{\begin{tabular}[t]{l}$\b x_{q(i)}$\end{tabular}}}}%
    \put(0.52244892,0.17869309){\makebox(0,0)[lt]{\lineheight{1.25}\smash{\begin{tabular}[t]{l}$d_i$\end{tabular}}}}%
    \put(0,0){\includegraphics[width=\unitlength,page=3]{images/curve_distance.pdf}}%
  \end{picture}%
\endgroup%
}
    \captionof{figure}{The notation introduced in section \ref{sec:reconstruction}.}
    \label{fig:distance_illustration}
  \end{Figure}
  
  Additionally, define
  $$ d_i = \metric(\b x_i, Y - U_i) $$
  where $\metric$ is the standard metric on $\R^2$.
  As long as
  \begin{equation}
    \label{eq:cont-distance}
    d_i \ge \max\{\metric(\b x_i, \b x_{p}), \metric(\b x_i, \b x_{q})\}
  \end{equation}
  and
  \begin{equation}
    \label{eq:discreet-distance}
    \argmin_{j = 0, \ldots, k-1} \metric(\b x_i, \b x_j) \in \{p, q\}
  \end{equation}
  for all $i = 0,\ldots, k$, we can find one of the neighboring points $\b x_{p}$ or $\b x_{q}$ of an arbitrary point $\b x_i$ by finding the nearest point among
  $ \b x_j $ for $j = 0 \ldots, i-1, i+1, \ldots, k-1$.
  
  When inequality~\eqref{eq:cont-distance} is not satisfied, the discretization alone does not provide sufficient information to determine which points of $X$ are in which part of the curve $\b\gamma$.
  A similar problem can occur if condition~\eqref{eq:discreet-distance} does not hold for some index, i.e. if the nearest neighbor to $\b x_i$ is not $\b x_{p}$ or $\b x_{q}$.
  \section{Algorithm and implementation}
  \label{sec:implementation}
  The surface reconstruction algorithm works in three steps. First we determine the correct order of given points, i.e.~we find a permutation $\sigma$ on $\{0, \ldots, k-1\}$, such that $ \sigma(i) < \sigma(j) $ implies that $t_i < t_j$ for all indices $0 \le i, j \le k-1$.
  In the second step, an approximation $\tilde{\b\gamma}$ of the starting curve $\b\gamma$ is obtained by fitting a cubic spline on ordered starting points.
  Finally, in the third step, we use the node positioning algorithm~\cite{slak2018generation} to discretize the curve for use in further calculations.
  
  \subsubsection*{Ordering the starting points}
  To find the appropriate permutation $\sigma$, firstly, the list of ordered points $\{\b x^\prime_i\}_{i=0}^{k-1}$ is initialized and an arbitrary starting point is chosen and assigned to the first position $\b x^\prime_0$.
  Now find the nearest neighbor $\b x_p$ of $\b x^\prime_0$, and assign it to $\b x^\prime_1$.
  The list is then build up inductively: Once $\b x^\prime_j$ is defined, find the nearest neighbor $\b x_p$ to $\b x^\prime_j$. If $\b x_p$ is not $\b x^\prime_{j-1}$, set $\b x^\prime_{j+1}$ to $\b x_p$. Otherwise let $\b x_p$ be the second nearest neighboring point. Now compare $\b x_p$ to $\b x^\prime_{j-2}$. If these are not equal, set $\b x^\prime_{j+1}$ to $\b x_p$. 
  The process is repeated until we cannot find a point $\b x_p$ that is not already in the ordered list.
  The process is also presented as pseudocode in Alg.~\ref{alg:enumeration}.
  
  In this paper, a $k$-d tree constructed from $X$ is used to query for nearest neighbors. Since it is more economical to store two arrays of indices rather than two arrays of vectors, the permutation is stored instead of a full list of ordered vectors when implementing the ordering procedure.
  
  \begin{algorithm}[H]
    \caption{Point enumeration}
    \SetAlgoLined
    \LinesNumbered
    \label{alg:enumeration}
    \KwData{An array of points on the plane $x[k]$.}
    \KwResult{An array representing a desired permutation $\sigma[k]$.}
    \KwResult{An array representing the inverse permutation $\sigma^{-1}[k]$.}
    \Begin {
    create integer arrays $\sigma[k]$, $\sigma^{-1}[k]$\;
    initialize k-d tree $T$ based on points $x$\;
    set $\sigma[0] \leftarrow 0$, $\sigma^{-1}[0] \leftarrow 0$\;
    \For{$j \leftarrow 0$ \KwTo $k-2$} {
    set $n_{\text{knn}} \leftarrow 2$\;
    \While{$\sigma[j+1]$ is not set}{
    find the $n_{\text{knn}}$-th nearest point $x[p]$ to $x[\sigma[j]]$ in $T$\nllabel{alg:enumeration:knn}\;
    \eIf{$p$ does not equal $\sigma[j-n_{\text{knn}}+1]$}{
    set $\sigma[j+1] \leftarrow p$, $\sigma^{-1}[p] \leftarrow j+1$\;
    }{
    increment $n_\text{knn}$ by one\;
    }
    }
    }
    }
  \end{algorithm}
  
  \subsubsection*{Domain shape reconstruction}
  The second step of the surface reconstruction fits a cubic spline $\tilde{\b \gamma}\colon \R \to \R^2$ to the points from the ordered list $\{\b x^\prime_i\}_i$. When $\tilde{\b \gamma}$ is obtained, all the surface shape information required is at our disposal - even between the discretization nodes.
  However to be able to reconstruct the entire domain $\Omega$ in the third and final step of the surface reconstruction algorithm, we must be able to distinguish between the interior and exterior of the curve $\tilde{\b\gamma}$.
  
  Suppose we are given a set of points $\{\b z_i\}_i \subset \partial\Omega$ and a set of accompanying normals $\{\b n_i\}_i \subset S^1$. To determine if an arbitrary point $\b z$ lies in $\Omega$, find the nearest point $\b z_i$ to $\b z$, and check if the vectors $\b z - \b z_i$
  and $\b n_i$ point in opposite directions. This can be done by computing the scalar product
  \begin{equation}
    \label{eq:scalar_product}
    \left\langle \b z - \b z_i, \b n_i \right\rangle.
  \end{equation}
  If the above scalar product~\eqref{eq:scalar_product} is negative, we conclude that $\b z \in \Omega$, otherwise not. 
  
  
  However, making any conclusions based on the sign of the equation~\eqref{eq:scalar_product} is not reliable. This basic idea typically fails in the proximity of sharp corners of $\b\gamma$ or more generally speaking, where the curve is not differentiable. We thus modify the algorithm to use the information provided by the interpolating spline $\tilde{\b\gamma}$.
  Let
  $$ s_0 < \cdots < s_{k-1} < s_k $$
  be the knots for the interpolated points~$\tilde{\b \gamma}(s_i) = \b x^\prime_i$ for all $i = 0, \ldots, k-1$ and $\tilde{\b \gamma}(s_k) = \b x^\prime_0$.
  For an arbitrary query point $\b x$, we find the closest point on the curve by minimizing the function
  \begin{equation}
    f(t) = \metric(\b x, \tilde{\b \gamma}(t)).
  \end{equation}
  Note that $f$ typically has many local minima. To obtain the correct one, we find the nearest point $\b x_q$ to $\b x$ among $X$ and use the inverse permutation $p = \sigma^{-1}(q)$. The desired value $t$ can now be approximated using bisection on the interval $[s_{p-1}, s_{p+1}]$. Suppose $t_\text{min}$ is the correct global minimum of $f$. Generally speaking $t_\text{min}$ is dependant on $\b x$, but for the sake of brevity we use $t_\text{min} = t_\text{min}(\b x)$ unless otherwise specified.
  
  Now the scalar product~\eqref{eq:scalar_product} is rewritten to take form 
  
  \begin{equation}
    \label{eq:int-problem}
    \left\langle \b x - \tilde{\b\gamma}(t_\text{min}), \tilde{\b \gamma}^{\prime\prime}(t_\text{min}) \right\rangle.
  \end{equation}
  Note that this procedure requires that the normals 'point outwards'.
  Our construction for $\tilde{\b \gamma}$ does not guarantee this, which we compensate for by introducing a constant $c$
  \begin{equation}
    c = -\sgn\left(\left\langle \b x_\text{int} - \tilde{\b\gamma}(t_\text{min}(\b x_\text{int})), \tilde{\b \gamma}^{\prime\prime}(t_\text{min}(\b x_\text{int})) \right\rangle\right),
  \end{equation}
  as the sign of the value of equation~\eqref{eq:int-problem} when applied to a point $\b x_\text{int}$ from the interior, i.e.\ $\b x_\text{int} \in \Omega$. Note the constant $c$ essentially flips the normals.
  
  The equation~\eqref{eq:int-problem} is then finally modified to
  \begin{equation}
    \label{eq:int-problem-const}
    \left\langle \b x - \tilde{\b\gamma}(t_\text{min}), c\tilde{\b \gamma}^{\prime\prime}(t_\text{min}) \right\rangle
  \end{equation}
  and enables us to determine if a point $\b x$ is in the interior of $\Omega$ or not.
  
  \subsubsection*{Discretization}
  The only remaining step for a complete surface reconstruction is to discretize the curve $\tilde{\b \gamma}$. The discretization is done by employing the node positioning algorithm proposed in~\cite{slak2018generation}. A detailed description of the node positioning algorithm used is out of the scope of this paper.

  \section{Problem setup}
  \label{sec:problem}
  We demonstrate the proposed surface reconstruction algorithm from chapter~\ref{sec:implementation} on a simplified moving-boundary problem -- a simulation of dendrite-like growth also known as \emph{solidification procedure}. A dendrite in metallurgy is a typical tree-like crystal structure that develops as molten metal solidifies~\cite{kobayashi1993modeling}. The dynamics of a real-life problem is thus mainly governed by the phase-transition physics from molten metal to a solid crystal structure.
  
  \subsection{Moving boundary}
  \label{sec:moving_boundary}
  Let the initial domain $\Omega = B_m / B_d$ be an annulus between a larger static circle $B_m$ with radius $R_m$ representing the boundary of the molten metal and smaller non-static circle $B_d$ initially with radius $R_d < R_m$ representing the dendrite's initial state
  \begin{align*}
    B_m & = \left \{ \b x \in \R^2, \ \left\| \b x \right\| \leq R_m \right \} \text{and} \\
    B_d & = \left \{ \b x \in \R^2, \ \left\| \b x \right\| \leq R_d \right \}.
  \end{align*}
  
  Generally speaking the velocity of the phase-transition front during the solidification is a function of the temperature field in the proximity. However, in this work, we simplify the problem to a degree, where this dependency is discarded -- the velocity is instead synthetically defined to result in a dendrite-like shape. To achieve that, all nodes $\b x_i$ from the $\partial B_d$ boundary are assigned a velocity $\b v_i = \b v_i(\b x_i)$ that depends on the position of the node $\b x_i$ and on the boundary normal $\b n_i = \b n_i(\b x_i)$
  \begin{equation}
    \label{eq:velocity}
    \b v_i = v_d\Big (\frac{1}{20} + \cos ^2 (2\phi_i)\Big )\b n_i.
  \end{equation}
  Here, $v_d$ is a constant model parameter and $\phi _i$ is the polar angle of $\b x_i$. This essentially means that at any given time step all the nodes from the dendrite boundary are moved $\b x _i ^{t+\mathrm{d}t} = \b x_i ^t + \mathrm{d}t\b v_i(\b x_i)$ which results in a dendrite-like growth.
  
  \subsection{Temperature field}
  \label{sec:temperature_field}
  Although the dynamics of the problem at hand is simplified by uncoupling the phase-transition front velocity and temperature field in the proximity, the latter is still computed at every time step, as shown in our implementation scheme in Fig.~\ref{fig:implementation_scheme}. 
  
  \begin{Figure}
    \centering
    \includegraphics[width=0.9\linewidth]{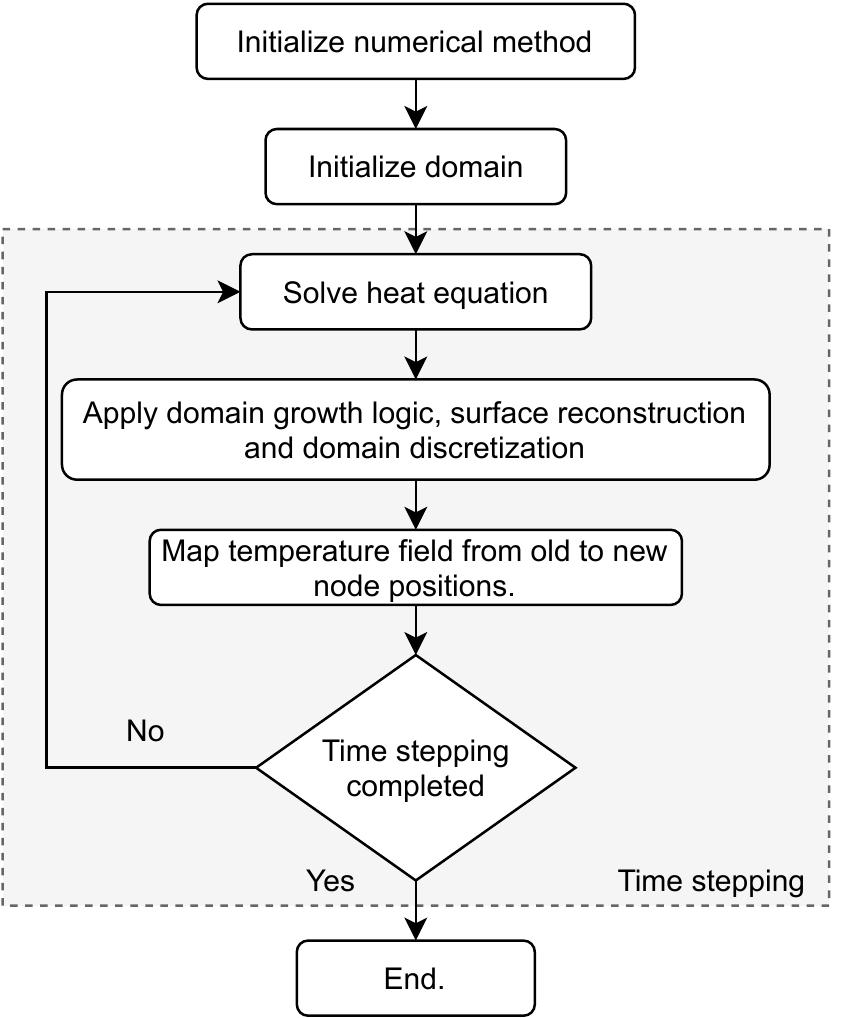}
    \captionof{figure}{Implementation scheme.}
    \label{fig:implementation_scheme}
  \end{Figure}
  
  The temperature field is governed by the dimensionless heat equation
  \begin{equation}
    \frac{\partial }{\partial t}T = \lap T
  \end{equation}
  The equation is discretized and finally takes the form
  \begin{equation}
    \label{eq:heat}
    T_{i + 1}(\b x) = T_i(\b x) + \mathrm{d}t \lap T_i(\b x)
  \end{equation}
  before it is numerically solved for all nodes $\b x \in \Omega$.
  
  \section{Results}
  \label{sec:results}
  The problem from section~\ref{sec:problem} is now simulated using our in-house Medusa library for meshless simulations. Since this is a theoretical problem, the simulation is done in a dimensionless sense.
  
  The outer radius $R_m$, representing the molten metal boundary, is constant and set to 1, while the initial radius of $B_d$ is set to $R_d=0.1$. On every time step all the dendrite boundary nodes are assigned a velocity as defined in equation~\eqref{eq:velocity}, where $v_d=0.04$. The temperature field is obtained on every time step before the domain growth logic is applied, as noted in implementation scheme in Fig.~\ref{fig:implementation_scheme}. After the surface had been reconstructed and nodes repositioned, it is important to map the temperature field from the old node positions to the new. This is achieved using the Inverse distance weighting (IDW), a procedure also known as Sheppard's interpolation~\cite{plestenjak2010numericne}.
  
  In this work, we used the RBF-FD with polyharmonic splines augmented with monomials of second order to compute the weights from equation~\eqref{eq:ansatz} and consequently compute the temperature field in the interior of $\Omega$. The temperatures at both boundaries, i.e.\ molten metal and dendrite boundary, are kept constant at 1 and 0 respectively.
  
  \begin{Figure}
    \centering
    \includegraphics[width=\linewidth]{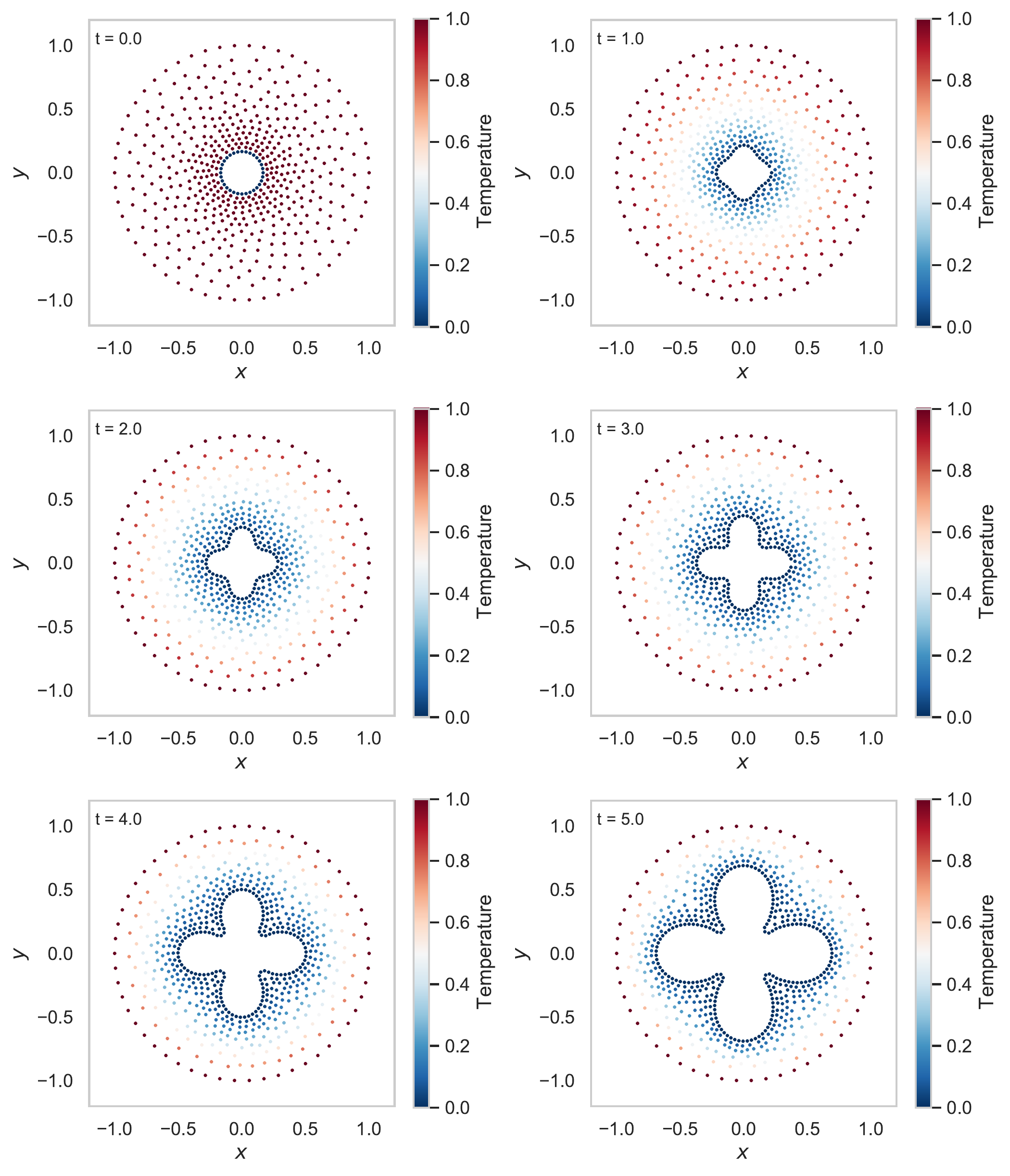}
    
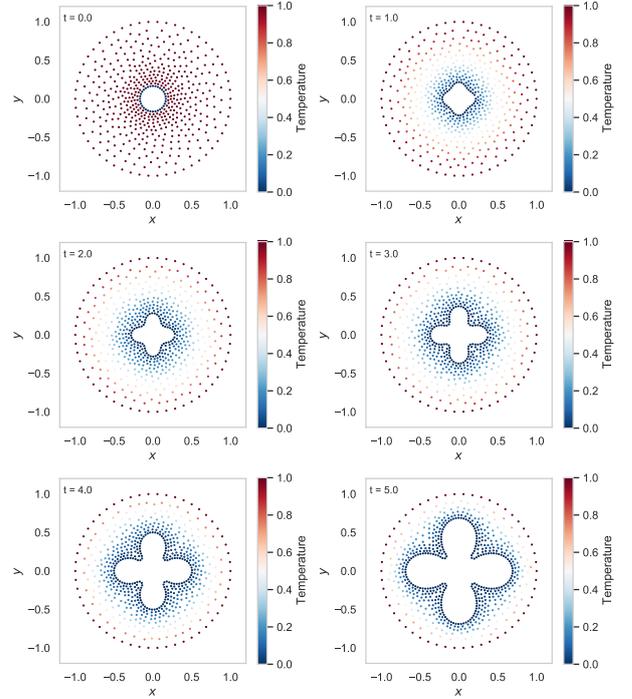
\captionof{figure}{Timelapse of dendrite-like growth using the proposed surface reconstruction algorithm.}
    \label{fig:animation}
  \end{Figure}
  
  The simulation is done by explicit time-marching for $N_t = 500$ time steps where a single time step is $\mathrm{d} t = 0.01$ units long, with a total simulation time $t_{tot} = 5$ units. The timelapse of the simulation is demonstrated in Fig.~\ref{fig:animation} for some selected simulation times. Also notice that the internodal distance $h(\b x)$ is smaller closer to the dendrite boundary, which is of great importance in a more realistic case where the dendrite tip velocity is a function of the temperature field in the proximity and thus needs to be accurately computed. In our case, the decreasing internodal distance $h$ provides us a with well defined and distinguishable dendrite shape.
  
  The number of nodes at the initialization time is $N=506$, while there are $N=637$ nodes at the simulation end. The number of nodes increases as the dendrite boundary is moving but also has a finer discretization compared to its surroundings. The total execution time is approximately 21 seconds on a \texttt{Windows Linux Subsystem (Ubuntu 20.04 LTS)} with \texttt{Intel(R) Core(TM) i7-9750H CPU @ 2.6GHz} and \texttt{16 GB} of RAM. The C++ code was compiled using \texttt{g++ (GCC) 9.3.0} for Linux with \texttt{-fopenmp -O3 -DNDEBUG} flags. 
  
  It is important to note that between two time steps all the domain shape information, apart from the discretization nodes, is discarded. Therefore, the domain shape is reconstructed at every single time step. In Fig.~\ref{fig:animation} we can observe how even after the 500 simulation time steps, the shape of the boundary remains smooth and symmetric. There are no visible irregularities or sights of unexpected discontinuities in the direction of normal vectors. This observation allows us to conclude that the proposed algorithm for surface reconstruction is stable.
  
  \section{Conclusions}
  \label{sec:conclusions}
  A detailed presentation of the proposed surface reconstruction algorithm is provided in this paper. In the first step ordering the discretized points in order is required, in the second step, cubic splines are used to reconstruct the surface shape and finally, in the third step, node positioning algorithm is used to obtain a new set of nodes on the boundary. We also explain how to distinguish between the interior and exterior of the reconstructed domain shape.
  
  Additionally, we demonstrate the proposed algorithm on a moving boundary problem, approximating the dendrite-like shape. At every time step, the surface is reconstructed and temperature field is computed using the RBF-FD approximation. However, in this paper, we discard the coupling between the moving dendrite tip velocity and the temperature field in the proximity. Instead, we assign a synthetically defined velocity to the boundary nodes simulating the dendrite-like growth. The next step could, therefore, be the removal of this simplification and implementation of the actual phase-transition physics making the simulation more realistic.
  
  In this paper, the surface reconstruction is executed at every time step, because we choose to discard all the available information about the domain shape between the discretization nodes after every time step. This is generally speaking unnecessary -- some domain information between two time steps doesn't change and could be simply propagated to the next time step.
  
  \section*{Acknowledgments}
  \label{sec:ack}
  The authors would like to acknowledge the financial
  support of the ARRS research core funding No.\ P2-0095.
  
  \bibliographystyle{unsrt}
  \bibliography{ref}
  
\end{multicols}

\end{document}